\documentclass[12pt,a4paper,reqno]{amsart}

\usepackage{amssymb,mathtools}

\edef\restoreparindent{\parindent=\the\parindent\relax}
\usepackage{parskip}
\restoreparindent

\makeatletter
\@namedef{subjclassname@2020}{\textup{2020} Mathematics Subject Classification}
\makeatother

\newtheorem{thm}{Theorem}
\newtheorem{lem}{Lemma}

\newtheorem{cor}{Corollary}
\newtheorem{conj}{Conjecture}

\newtheorem{Thm}{Theorem}

\theoremstyle{definition}
\newtheorem{rem}{Remark}
\newtheorem{defn}{Definition}

\newcommand{\IC}{{\mathbb C}}
\newcommand{\ID}{{\mathbb D}}

\newcommand{\IT}{{\mathbb T}}
\newcommand{\real}{{\operatorname{Re}\,}}
\newcommand{\ima}{{\operatorname{Im}\,}}
\newcommand{\eit}{{e^{i\theta}}}
\newcommand{\hl}{{{\bf h} \log^+{\bf h}}}

\newcommand{\be}{\begin{equation}}
\newcommand{\ee}{\end{equation}}

\newcommand{\blem}{\begin{lem}}
\newcommand{\elem}{\end{lem}}

\newcommand{\bdefn}{\begin{defn}}
\newcommand{\edefn}{\end{defn}}

\newcommand{\bthm}{\begin{thm}}
\newcommand{\ethm}{\end{thm}}

\newcommand{\bcor}{\begin{cor}}
\newcommand{\ecor}{\end{cor}}

\newcommand{\bconj}{\begin{conj}}
\newcommand{\econj}{\end{conj}}

\newcommand{\brem}{\begin{rem}}
\newcommand{\erem}{\end{rem}}

\newcommand{\bpf}{\begin{proof}}
\newcommand{\epf}{\end{proof}}

\begin{document}
	
\bibliographystyle{abbrv}

\title[Zygmund's theorem for harmonic quasiregular mappings]
{Zygmund's theorem for harmonic quasiregular mappings}

\author{Suman Das}
\address{Suman Das\vskip0.05cm Department of Mathematics with Computer Science, Guangdong Technion - Israel
	Institute of Technology, Shantou, Guangdong 515063, P. R. China.}
\email{suman.das@gtiit.edu.cn}

\author{Jie Huang}
\address{Jie Huang \vskip0.05cm Department of Mathematics with Computer Science, Guangdong Technion - Israel
	Institute of Technology, Shantou, Guangdong 515063, P. R. China.}
\email{jie.huang@gtiit.edu.cn}

\author{Antti Rasila}
\address{Antti Rasila \vskip0.05cm Department of Mathematics with Computer Science, Guangdong Technion - Israel
	Institute of Technology, Shantou, Guangdong 515063, P. R. China. \vskip0.025cm Department of Mathematics, Technion - Israel
	Institute of Technology, Haifa 3200003, Israel.}
\email{antti.rasila@gtiit.edu.cn; antti.rasila@iki.fi}



\begin{abstract}
Given an analytic function $f=u+iv$ in the unit disk $\ID$, Zygmund's theorem gives the minimal growth restriction on $u$ which ensures that $v$ is in the Hardy space $h^1$. This need not be true if $f$ is a complex-valued harmonic function. However, we prove that Zygmund's theorem holds if $f$ is a harmonic $K$-quasiregular mapping in $\ID$. Our work makes further progress on the recent Riesz-type theorem of Liu and Zhu (Adv. Math., 2023), and the Kolmogorov-type theorem of Kalaj (J. Math. Anal. Appl., 2025), for harmonic quasiregular mappings. We also obtain a partial converse, thus showing that the proposed growth condition is the best possible. Furthermore, as an application of the classical conjugate function theorems, we establish a harmonic analogue of a well-known result of Hardy and Littlewood.
\end{abstract}

\maketitle
\pagestyle{myheadings}
\markboth{S. Das, J. Huang, and A. Rasila}{Zygmund's theorem for harmonic quasiregular mappings}



\section{Introduction}\label{sec1}

\subsection{Notations and preliminaries}
Suppose $\ID\coloneqq\{z \in {\mathbb C}:\, |z|<1\}$ is the open unit disk in the complex plane $\IC$, and $\mathbb{T}\coloneqq\{z \in \IC : |z|=1\}$ is the unit circle. For a function $f$ analytic in $\ID$,
the \textit{integral means} of $f$ are defined as $$M_p(r, f) \coloneqq \left( \frac{1}{2\pi}\int_{0}^{2\pi}\vert f(r e^{i\theta}) \vert ^p d\theta \right)^{1/p} \quad \text{for} \quad 0 < p < \infty,$$ and $$M_\infty(r,f) \coloneqq \sup_{ \vert z \vert=r} \vert f(z) \vert.$$
The function $f$ is said to be in the \textit{Hardy space} $H^p$ $(0 < p \le \infty)$ if $$\|f\|_p \coloneqq \lim_{r \to 1^-} M_p(r,f)<\infty.$$ Every function $f \in H^p$ has ``nice" boundary behaviour, in the sense that the radial limit $$f(e^{i\theta})\coloneqq \lim_{r \to 1^-} f(r\eit)$$ exists in almost every direction, and is of class $L^p\coloneqq L^p(\IT)$. Detailed surveys on Hardy spaces can be found, e.g., in the books of Duren \cite{Duren} and Koosis \cite{Koosis}. Throughout this paper, we follow notations from \cite{Duren}.

A real-valued function $u(x,y)$, twice continuously differentiable in $\ID$, is called \textit{harmonic} if it satisfies the Laplace equation $$\Delta u =\frac{\partial^2 u}{\partial x^2}+\frac{\partial^2 u}{\partial y^2}= 0 \quad \text{in} \quad \ID.$$ For a complex-valued function $f=u+iv$, the complex partial derivatives have the form $$f_z=\frac{1}{2} \left(f_x-if_y\right) \quad \text{and} \quad f_{\bar{z}}=\frac{1}{2} \left(f_x+if_y\right),$$ where $z=x+iy\in \IC$. In view of this, it is easy to see that the Laplacian operator $\Delta$ can also be written as $$\Delta=4\, \frac{\partial^2}{\partial \bar{z} \partial z} \cdot$$

A function $f=u+iv$ is harmonic in $\ID$, if $u$ and $v$ are real-valued harmonic functions in $\ID$. Every such function has a representation $f=h+\overline{g}$, where $h$ and $g$ are analytic in $\mathbb{D}$. This representation is unique up to an additive constant. In this paper we usually assume $g(0)=0$, unless mentioned otherwise. Similar to the $H^p$ spaces, the harmonic Hardy spaces $h^p$ are defined as the class of functions $f$, harmonic in $\ID$, such that $$\|f\|_p=\lim_{r\to 1^-} M_p(r,f) < \infty.$$ We refer to the books of Duren \cite{Duren:Harmonic} for the theory of planar harmonic mappings and Pavlovi\'c \cite{Pavbook} for a concise survey on $h^p$ spaces.

Given $K\ge 1$, a sense-preserving harmonic function $f=h+\bar{g}$ is said to be $K$-\textit{quasiregular} if its \textit{complex dilatation} $\omega = g'/h'$ satisfies the inequality $$|\omega(z)| \le k <1 \quad (z\in \ID),$$ where $$k\coloneqq\frac{K-1}{K+1}.$$ We say that $f$ is $K$-\textit{quasiconformal} if $f$ is $K$-quasiregular and homeomorphic in $\ID$. One can find the $H^p$-theory for quasiregular and quasiconformal mappings, for example, in the survey of Astala and Koskela \cite{ast_kos} and the paper of Adamowicz and Gonz\'alez \cite{Adam_Gonz}. Some recent results on the Hardy space of harmonic quasiconformal mappings are obtained in \cite{wang_rasila}.

\subsection{Classical results of Riesz, Kolmogorov and Zygmund}
Given a real-valued harmonic function $u$ in $\ID$, let $v$ be its harmonic conjugate, normalized by the condition $v(0)=0$. It is a natural question that if $u$ has a certain property, whether so does $v$. In the context of boundary behaviour, M.~Riesz established this principle in a remarkably precise form:
\begin{Thm}\label{MRiesz}\cite[Theorem 4.1]{Duren}
If $u \in h^p$ for some $p$, $1<p<\infty$, then its harmonic conjugate $v$ is also of class $h^p$. Furthermore, there is a constant $A_p$, depending only on $p$, such that $$M_p(r,v) \le A_p \, M_p(r,u),$$ for all $u \in h^p$.
\end{Thm}
Curiously, the theorem fails for the cases $p=1$ and $p=\infty$, examples are given in \cite[p. 56]{Duren}. Although the harmonic conjugate of an $h^1$ function need not be in $h^1$, Kolmogorov crucially pointed out that it does belong to $h^p$ for all $p<1$. The result, in a stronger form, is as follows.
\begin{Thm}\label{Kolmo}\cite[Theorem 4.2]{Duren}
If $u \in h^1$, then its harmonic conjugate $v\in h^p$ for all $p<1$. Furthermore, there is a constant $B_p$, depending only on $p$, such that $$M_p(r,v) \le B_p \, M_1(r,u),$$ for all $u \in h^1$.
\end{Thm}
It follows from the theorems of Riesz and Kolmogorov that the condition $u\in h^1$ is not enough to ensure $v\in h^1$, but the slightly stronger hypothesis $u\in h^p$ (for some $p>1$) is sufficient. Therefore, one naturally asks for the ``minimal" growth restriction on $u$ which should imply $v\in h^1$. As established by Zygmund, such a condition is the boundedness of the integrals $$\int_{0}^{2\pi} |u(r\eit)| \log^+ |u(r\eit)| \, d\theta,$$ where $\log^+ x=\max\{\log x,\, 0\}$. The class of such harmonic functions $u$ in $\ID$ is denoted by $\hl$. The following is the precise statement of Zygmund's theorem.
\begin{Thm}\label{Zyg}\cite[Theorem 4.3]{Duren}
If $u \in \hl$, then its conjugate $v$ is of class $h^1$, and
$$M_1(r,v) \le \frac{1}{2\pi}\int_{0}^{2\pi} |u(r\eit)| \log^+ |u(r\eit)| \, d\theta + 3e.$$
\end{Thm}
It should be noted that $h^p \subset \hl$ for every $p>1$, so that Zygmund's theorem is indeed best possible. We refer the reader to the paper of Pichorides \cite{Pichorides} for the sharp constants in the Riesz, Kolmogorov, and Zygmund theorems.

\subsection{Recent developments for harmonic quasiregular mappings}
Suppose $f=u+iv$ is a harmonic function in $\ID$. If $u\in h^p$, $p>1$, then the imaginary part $v$ does not necessarily belong to the space $h^p$, i.e., Riesz theorem is not true for harmonic functions. The question arises: Under what additional condition(s) does a harmonic analogue of the Riesz theorem hold? Recently, Liu and Zhu established in their breakthrough paper \cite{Liu_Zhu} that such a condition is the quasiregularity of $f$.
\begin{Thm}\label{LZ}\cite{Liu_Zhu}
Let $f=u+iv$ be a harmonic $K$-quasiregular mapping in $\ID$ such that $u \ge 0$ and $v(0)=0$. If $u \in h^p$ for some $p\in (1,2]$, then also $v$ is in $h^p$. Furthermore, there is a constant $C(K,p)$, depending only on $K$ and $p$, such that $$M_p(r,v) \le C(K,p)M_p(r,u).$$ Moreover, if $K=1$, i.e., $f$ is analytic, then $C(1,p)$ coincides with the optimal constant in the Riesz theorem.
\end{Thm}
It is important to note that the assumption ``$f$ is harmonic in $\ID$" cannot be removed. There exists a quasiconformal mapping for which the first coordinate function belongs to the Hardy space, but the whole function does not (see \cite[p. 36]{ast_kos}).

Later, Kalaj \cite{Kalaj_Kolmogorov} refined Theorem \ref{LZ} by removing the condition $u\ge 0$. Furthermore, he proved a Kolmogorov type theorem for harmonic quasiregular mappings in $\ID$, thereby making considerable progress in this direction.
\begin{Thm}\label{Kal_Kol}\cite{Kalaj_Kolmogorov}
Suppose $f=u+iv$ is a harmonic $K$-quasiregular mapping in $\ID$ such that $u>0$ and $v(0)=0$. Then $u\in h^1$ and $v \in h^p$ for every $p<1$. Furthermore, they satisfy the inequalities
$$M_p^p(r,v) \le \sec \left(\frac{p\pi}{2}\right) \left(K^2 M_1^p(r,u)-(K^2-1)M_p^p(r,u)\right)$$ and $$(2-K^2)M_1^p(r,u)\le (2-K^2)M_p^p(r,u)+\cos\left(\frac{p\pi}{2}\right) M_p^p(r,v).$$ All the constants are asymptotically sharp as $K\to 1$.
\end{Thm}
Kalaj also points out that the second inequality is new even for analytic functions. Several other Riesz type theorems for harmonic functions can be found in the recent works of Kalaj \cite{Kalaj, kalaj_conj}, Chen and Hamada \cite{ChenHam}, Melentijevi\'c \cite{Mel1}, and Melentijevi\'c and Markovi\'c \cite{Mel2}. Furthermore, Chen and Huang \cite{chen_huang} have proved Riesz theorem for pluriharmonic functions where the assumption $u\ge 0$ is dropped as well.

\subsection{Main results}
In this paper, we make further advance along this line and establish a Zygmund type theorem for harmonic quasiregular mappings. For a harmonic $K$-quasiregular mapping $f=u+iv$ in $\ID$, we find out the growth condition on $u$ that would imply $f\in h^1$ (and $v\in h^1$). The following is the central result of this paper.

\bthm\label{DHR1}
Let $f=u+iv$ be a harmonic $K$-quasiregular mapping in $\ID$ such that $u\ge 1$ (or $u\le -1$) and $v(0)=0$. If $u \in \hl$, then $f \in h^1$ and
\be\label{DHR_eq1}
M_1(r,f) \le \frac{K^2}{2\pi} \int_{0}^{2\pi} |u(r\eit)| \log^+ |u(r\eit)| \, d\theta + |u(0)|\left[1-K^2 \log|u(0)|\right].
\ee
Consequently, $v \in h^1$ and $M_1(r,v)$ satisfies inequality \eqref{DHR1}.
\ethm

\brem
The conclusion $f\in h^1$ is true under the slightly weaker hypothesis that $u\ge C$ for some constant $C$, since adding a constant to $u$ does not affect $v$ or the quasiregularity of $f$, although we do not know if the estimate \eqref{DHR_eq1} still holds in this case. However, we suspect that the assumption $u \ge 1$ (or $u \ge C$) may actually be redundant, but a proof remains elusive.
\erem

Theorem \ref{DHR1} is the best possible, in the sense that the growth restriction imposed on $u$ cannot be weakened. This can be seen from the following partial converse, which surprisingly does not require the quasiregularity hypothesis.

\bthm\label{DHR2}
Let $f=h+\bar{g}$ be a harmonic function in $\ID$ such that $\real f\ge C$ for some constant $C$, and $\ima h$ is non-vanishing. If $f \in h^1$, then $\real f \in \hl$.
\ethm
The proofs of the main theorems are presented in Section \ref{sec2}. In Section \ref{sec3}, we study the harmonic analogue of a well-known result of Hardy and Littlewood as an application of the classical Riesz and Kolmogorov theorems.

\section{Proof of the main theorems}\label{sec2}


\subsection{Proof of Theorem \ref{DHR1}}
Our technique is based on a method of Green's theorem, which has also been used in \cite{Kalaj_Kolmogorov} and \cite{Liu_Zhu}.

First assume that $u \ge 1$, so that $\log^+|u(r\eit)|=\log u(r\eit)$. Let us write $f=h+\bar{g}$, and let $\omega(z)=g'(z)/h'(z)$ be its complex dilatation. An elementary computation gives \begin{align*}\Delta\{u\log u\} & = \left(\frac{\partial^2}{\partial x^2}+\frac{\partial^2}{\partial y^2}\right)\{u \log u\} = \frac{u_x^2+u_y^2}{u}.\end{align*}
Suppose \(F=h+g\). Then \[\real F= \real h + \real g = \real f =u,\] so that \[F'=h'+g'=u_x-i\, u_y.\] This implies \[u_x^2+u_y^2=|h'+g'|^2.\]
Therefore, we have \begin{align*}
\Delta\{u\log u\} = \frac{|h'|^2\,|1+\omega|^2}{u} \ge (1-k)^2\, \frac{|h'|^2}{u},
\end{align*} as $|\omega(z)| \le k$. On the other hand, one finds \begin{align*}\Delta\{|f|\} & = 4\, \frac{\partial^2}{\partial \bar{z} \partial z} \{|f|\}\\ & = 2\, \frac{\partial}{\partial \bar{z}} \left(\frac{h'(\bar{h}+g)+(h+\bar{g})g'}{|f|}\right)\\ & = 2\, \frac{|f|(|h'|^2+|g'|^2)-\left(h'(\bar{h}+g)+(h+\bar{g})g'\right)\frac{\overline{g'}(\bar{h}+g)+(h+\bar{g})\overline{h'}}{2|f|}}{|f|^2}\\&=\frac{2\,|f|^2(|h'|^2+|g'|^2)-\left(|f|^2(|h'|^2+|g'|^2)+2\real (\overline{h'}g'f^2)\right)}{|f|^3}\\& = \frac{|f|^2(|h'|^2+|g'|^2) - 2 \real (\overline{h'}g'f^2)}{|f|^3}.\end{align*} 
Thus, it follows that \begin{align*}\Delta\{|f|\} & \le \frac{|f|^2(|h'|^2+|g'|^2) + 2 |h'||g'||f|^2}{|f|^3}\\ & = \frac{(|h'|+|g'|)^2}{|f|} = (1+|\omega|)^2\, \frac{|h'|^2}{|f|}\\ & \le (1+k)^2\, \frac{|h'|^2}{|f|}. \end{align*} These inequalities imply \begin{align*}
\frac{\Delta\{|f|\}}{\Delta\{u\log u\}} \le \left(\frac{1+k}{1-k}\right)^2 \frac{u}{|f|} \le \left(\frac{1+k}{1-k}\right)^2 = K^2,
\end{align*}
and therefore, \be \label{DHR_eq2}\Delta\{|f|\} \le K^2 \Delta\{u\log u\}.\ee The next step is to apply Green's theorem in the form $$r \int_{0}^{2\pi} \frac{\partial \varphi}{\partial r} \, d\theta = \iint \limits_{|z|\le r} \Delta \varphi \, dx\,dy,$$ where $z=x+iy$. Then, it follows from \eqref{DHR2} that $$\frac{d}{dr} \int_{0}^{2\pi} |f(r\eit)|\, d\theta \le K^2 \frac{d}{dr} \int_{0}^{2\pi} u(r\eit) \log u(r\eit)\, d\theta.$$ Now, integrating from $0$ to $r$ and recalling that $v(0)=0$, we get $$\int_{0}^{2\pi} |f(r\eit)|\, d\theta - 2\pi u(0) \le K^2 \left[\int_{0}^{2\pi} u(r\eit) \log u(r\eit)\, d\theta - 2\pi u(0) \log u(0)\right].$$ Therefore $$M_1(r,f) \le \frac{K^2}{2\pi} \int_{0}^{2\pi} u(r\eit) \log u(r\eit) \, d\theta + u(0)\left[1-K^2 \log u(0)\right],$$ which is the desired inequality.

For the case $u \le -1$, let us write $F=e^{i\pi}f$. We see that $F$ is harmonic $K$-quasiregular in $\ID$ with $u_F = \real F = -u \ge 1$, and $v_F = \ima F = -v$ satisfying $v_F(0)=0$. Then, from what we have already proved, it follows that $$M_1(r,F) \le \frac{K^2}{2\pi} \int_{0}^{2\pi} u_F(r\eit) \log u_F(r\eit) \, d\theta + u_F(0)\left[1-K^2 \log u_F(0)\right].$$ Since $M_1(r,f) = M_1(r,F)$ and $u_F=|u|$, the proof is complete.\qed

The proof of the converse (Theorem \ref{DHR2}) relies on the following lemma, which is of some independent interest.

\blem\label{lem1_DHR}
Let $f=h+\bar{g}$ be a harmonic function in $\ID$ such that $\ima h$ is non-vanishing. If $f\in h^1$, then the analytic function $F=h+g$ is in $H^1$.
\elem

\bpf
Suppose $f \in h^1$, and assume that $h=a+ib$ and $g=c+id$. Then $$\real f = \real F = a+c\, \in h^1,$$ and $$\ima f = b-d\, \in h^1,$$ since $\real f,\, \ima f \le |f|$. We recall \cite[Theorem 1.1]{Duren} that every real-valued harmonic function in $h^1$ can be expressed as the difference of two positive harmonic functions. Let us write $$b-d = u_1 - u_2,$$ where $u_1$ and $u_2$ are positive harmonic functions in $\ID$. Since $\ima h = b$ is continuous and non-vanishing, either $b>0$ or $b<0$ in the whole disk $\ID$. If $b>0$, then $$d=(b+u_2)-u_1,$$ which is again the difference of two positive harmonic functions, and hence is in $h^1$. Similarly, if $b<0$, we have $$d=u_2-(u_1-b) \in h^1.$$ In either case, $$\ima F = b+d = 2d + (u_1-u_2) \in h^1.$$ Therefore, $F \in H^1$. This completes the proof.
\epf

\brem
\begin{enumerate}
	\item[(i)] The usual normalisation $g(0)=0$ implies $\ima g (0)=0$. However, if we drop the condition $g(0)=0$, then Lemma \ref{lem1_DHR} is true if either $\ima h$ or $\ima g$ is non-vanishing, and the proof is essentially the same.
	\item[(ii)] If $f \in h^1$, it seems that the corresponding $F$ need not be in $H^1$. In view of this, the hypothesis $\ima h \neq 0$ (or $\ima g \neq 0$) seems to be a surprisingly simple assumption that ensures $F \in H^1$. A similar condition appears in Kalaj's refinement of the Liu-Zhu result (Theorem 2.2 of \cite{Kalaj_Kolmogorov}), where it is assumed that $f(\ID) \cap (-\infty,0) = \emptyset$. 
\end{enumerate}
\erem

\subsection{Proof of Theorem \ref{DHR2}}
Suppose $f=u+iv$. Without any loss of generality, we may assume that $u \ge 1$, as adding a constant to $u$ does not affect $v$ or the condition $f\in h^1$. Since $\ima h$ is non-vanishing, by Lemma \ref{lem1_DHR} we have $F=h+g \in H^1$. Let us write $$F(z)=Re^{i\Phi} \quad \left(|\Phi| \le \frac{\pi}{2}\right),$$ where $z=r\eit$, $R=R(r,\theta)$ and $\Phi=\Phi(r,\theta)$. Then $$\real F =\real f = u=R \cos (\Phi).$$ Now, the function $F(z)\log F(z)$ is analytic in $\ID$, with $$\real[F(z)\log F(z)] = R \cos (\Phi) \log (R) - \Phi R \sin (\Phi).$$ We apply mean value theorem to this real part to obtain \begin{align*}
\int_{0}^{2\pi} \left[R \cos (\Phi) \log (R) - \Phi R \sin (\Phi)\right]\, d\theta & = 2\pi \real \left[F(0) \log F(0)\right]\\ & = 2\pi \real \left[h(0) \log h(0)\right],
\end{align*}
since $g(0)=0$. It is to be noted that $\ima h \neq 0$ implies $h \neq 0$, so that $\log h(0)$ is well-defined. Thus we have $$\int_{0}^{2\pi} R \cos (\Phi) \log (R) \, d\theta = \int_{0}^{2\pi} \Phi R \sin (\Phi)\, d\theta + 2\pi \real \left[h(0) \log h(0)\right].$$ It follows that \begin{align*}
\int_{0}^{2\pi} u(r\eit) \log u(r\eit)\, d\theta & \le \int_{0}^{2\pi} R\cos (\Phi) \log (R)\, d\theta\\ & \le |\Phi| \int_{0}^{2\pi} R \, d\theta + 2\pi \left| \real \left[h(0) \log h(0)\right]\right|\\ & \le \frac{\pi}{2} \int_{0}^{2\pi} |F(r\eit)|\, d\theta +2\pi |h(0)\log h(0)|, \end{align*} which is finite since $F \in H^1$. Therefore, $u \in \hl$ and the proof is complete.\qed

\brem
In Theorem \ref{DHR1}, we assume that $v(0) =0$, while in Theorem \ref{DHR2}, the condition $\ima h \neq 0$ implies $v(0) \neq 0$ (as $g(0)=0$). This is not really an inconsistency, since adding a constant to $v$ (therefore, to $f$) does not affect $u$ or the condition $f\in h^1$. So Theorems \ref{DHR1} and \ref{DHR2} are simultaneously valid.
\erem

\section{An application of the theorems of Riesz and Kolmogorov}\label{sec3}
Here we discuss how the classical theorems of Riesz and Kolmogorov can be applied to obtain numerous results for planar harmonic Hardy spaces. We give an example of a famous theorem of Hardy and Littlewood, that seems to be unknown for complex-valued harmonic functions.

Let $f$ be an analytic function in $\ID$ with $f(z)=\sum a_n z^n$. If $f \in H^p$ for some $p$, it is interesting to understand the behaviour of the Taylor coefficients $a_n$. We should note that it makes sense only to talk about the eventual behaviour as $n\to \infty$, since any finite number of coefficients can be changed arbitrarily without affecting the fact $f\in H^p$. The following theorem of Hardy and Littlewood describes the Taylor coefficients of an $H^p$-function in a precise form.

\begin{Thm}\label{HL_Taylor}\cite[Theorem 6.2]{Duren}
Let $0<p\le 2$ and $$f(z)=\sum_{n=0}^{\infty} a_n\, z^n \in H^p.$$ Then we have $$\sum_{n=0}^{\infty} (n+1)^{p-2} |a_n|^p < \infty.$$
\end{Thm}
The converse of this theorem is false for $0<p<2$ but holds for the indices $p\ge 2$, and can be obtained from Theorem \ref{HL_Taylor} via a duality argument. However, Pavlovi\'c \cite{Pav_New} showed that if the sequence $\{a_n\}$ decreases to $0$, then the converse is indeed true for all $p\ge 1$.

\begin{Thm}\label{HL_conv}\cite{Pav_New}
Suppose $1\le p<\infty$. If $a_n \downarrow 0$ as $n\to \infty$, then the function $f(z)=\sum a_n z^n$ is in $H^p$ if and only if $$\sum_{n=0}^{\infty} (n+1)^{p-2} a_n^p < \infty.$$
\end{Thm}

Similar problem for planar harmonic functions is not known, but can be easily deduced from the conjugate function theorems. It is worth mentioning that the case $p \le 1$ is usually difficult to handle for harmonic functions. Nevertheless, we can produce a slightly different result that deals with this case (part (ii) of Theorem \ref{DHR3}).

Let $f=h+\bar{g}$ be a harmonic function in $\ID$ with \be \label{DHR_eq3}h(z)=\sum_{n=0}^{\infty} a_n\, z^n, \quad \text{and} \quad g(z)=\sum_{n=0}^{\infty} b_n\, z^n.\ee  We should note that $b_0=0$ as $g(0)=0$, but for our forthcoming discussion, it is convenient to write the sum as $\sum_{n=0}^\infty$.

\bthm\label{DHR3}
If $f$ is a harmonic functions in $\ID$ with the series representation \eqref{DHR_eq3}, then the following statements hold:
\begin{enumerate}
	\item[(i)] If $f\in h^p$ for some $p\in (1,2]$, then $$\sum_{n=0}^{\infty} (n+1)^{p-2} \left(|a_n|^p+|b_n|^p\right) < \infty.$$ Furthermore, suppose $a_n\downarrow0$ and $b_n\downarrow0$. Then $f\in h^p$ $(p>1)$ if and only if $$\sum_{n=0}^{\infty} (n+1)^{p-2} \left(a_n^p+b_n^p\right) < \infty.$$
	
	\item[(ii)] If $f \in h^1$, then $$\sum_{n=0}^{\infty} (n+1)^{p-2} \left(|a_n|^p+|b_n|^p\right) < \infty$$ for every $p<1$.
	
	\item[(iii)] Suppose $\ima h \neq 0$ and $a_n\downarrow0$, $b_n\downarrow0$. If $f\in h^1$, then $$\sum_{n=0}^\infty \frac{a_n+b_n}{n+1} < \infty.$$
\end{enumerate}
\ethm

\bpf
(i) The main step is to show that $h, g \in H^p$. This is contained in Lemma 2.2 of \cite{Liu_Zhu}, but we give a slightly different proof based on ideas from the more general result \cite[Lemma 1]{DK2}. Let us write $f=u+iv$, where $u$ and $v$ are real-valued harmonic functions in $\ID$. Suppose $u_1$ is the harmonic conjugate of $u$ normalised by $u_1(0)=v(0)$, while $v_1$ is the harmonic conjugate of $v$ with $v_1(0)=-u(0)$. Assume that $U=u+iu_1$ and $V=v+iv_1$. Then $U$, $V$ are analytic in $\ID$, and $$f=\real U + i\,\real V = \frac{1}{2} \left(U+\overline{U}\right)+\frac{i}{2} \left( V+\overline{V}\right)=\frac{1}{2}\left(U+iV\right)+\frac{1}{2}\overline{\left(U-iV\right)}.$$ Moreover, the above normalisations suggest that $$\left(\frac{U-iV}{2}\right)(0)=0,$$ so that the uniqueness of $h$ and $g$ implies $$h=\frac{1}{2}(U+iV)\quad \text{and} \quad g=\frac{1}{2}(U-iV).$$ As $\vert u\vert, \vert v\vert \le \vert f\vert$, clearly $u,v\in h^p$. In view of Riesz theorem (Theorem \ref{MRiesz}), we have $u_1,v_1 \in h^p$. It follows that $U$, $V$ are in $H^p$, therefore so are $h$ and $g$.

Now, we apply Theorem \ref{HL_Taylor} to the functions $h$ and $g$ to find $$\sum_{n=0}^{\infty} (n+1)^{p-2} |a_n|^p < \infty \quad \text{and} \quad \sum_{n=0}^{\infty} (n+1)^{p-2} |b_n|^p < \infty,$$ for $1<p\le 2$. Combining these two inequalities, we get the desired result.

The proof for the second part of the theorem is identical, except appealing to Theorem \ref{HL_conv}. The details are omitted.

(ii) The reasoning is similar to that in part (i), except for the fact that we need to use Kolmogorov theorem (Theorem \ref{Kolmo}), instead of Riesz theorem. If $f\in h^1$, using an analogous argument as above we find that $h,g$ are in $H^p$ for every $p<1$. Therefore, Theorem \ref{HL_Taylor} gives $$\sum_{n=0}^{\infty} (n+1)^{p-2} \left(|a_n|^p+|b_n|^p\right) < \infty,$$ for $0<p<1$.

(iii) If $f \in h^1$ and $\ima h$ is non-vanishing, Lemma \ref{lem1_DHR} shows that the analytic function $F\in H^1$, where $$F(z) =h(z)+g(z)= \sum_{n=0}^{\infty} (a_n+b_n)\, z^n.$$ If $a_n$ and $b_n$ decrease to $0$, so does $a_n+b_n$. Therefore, it follows from Theorem \ref{HL_conv}, for $p=1$, that $$\sum_{n=0}^\infty \frac{a_n+b_n}{n+1} < \infty.$$ The proof of the theorem is complete.
\epf

\brem
It is obvious that the converse of Theorem \ref{DHR3}\,(iii) holds even without the condition $\ima h \neq 0$. Indeed, if $$\sum_{n=0}^\infty \frac{a_n+b_n}{n+1} < \infty,$$ one sees from Theorem \ref{HL_conv} that $h,\,g\in H^1$. Thus, $f \in h^1$.
\erem

\subsection*{Acknowledgements} The first-named author thanks Petar~Melentijevi\'c for helpful discussions on the space $h^1$ that led to the classical theory of Zygmund. The research was partially supported by the Natural Science Foundation of Guangdong Province (Grant no. 2024A1515010467) and the Li~Ka~Shing Foundation STU-GTIIT Joint Research Grant (Grant no. 2024LKSFG06).


\bibliography{references}


\end{document}